\documentclass[12pt]{amsart}
\usepackage{amssymb}
\theoremstyle{plain}
\newtheorem
{thm}{Theorem}[section]
\newtheorem
{proposition}[thm]{Proposition}
\newtheorem
{lemma}[thm]{Lemma}

\newtheorem
{corollary}[thm]{Corollary}

\newcommand{\ka}{\kappa}\theoremstyle{definition}

\newcommand{\la}{\lambda}
\newcommand{\gm}{\gamma}
\newcommand{\Gm}{\Gamma}
\newcommand{\thi}{\vartheta_{1}}
\newcommand{\tht}[2]{\theta_{#1,#2}^{s}}

\newcommand{\Ref}[1]{{(\ref{#1})}}

\newcommand{\half}{\frac{1}{2}}
\newcommand{\lt}{(\la,\tau)}

\renewcommand{\Im}{\operatorname{Im}}
\newcommand{\Z}{\mathbb{Z}}
\theoremstyle{definition}

\title[Elliptic Selberg Integrals and Conformal Blocks]
{Elliptic Selberg Integrals and Conformal Blocks}
\author[G. Felder, L. Stevens, and A. Varchenko]
{G. Felder $^{\,\star}$,
L. Stevens$^{\, \diamond}$, and
A. Varchenko$^{\,\diamond,1}$}
\thanks{$^1$ Supported in part by NSF grant DMS-9801582}
\begin{document}
\maketitle
\begin{center}
{\it
$^\star$ Departement Mathematik, ETH-Zentrum, 8092 Z\"urich, Switzerland,

felder@math.ethz.ch

\medskip

$^\diamond$ Department of Mathematics, University of North Carolina
at Chapel Hill,

Chapel Hill, NC 27599-3250, USA,

lstevens@email.unc.edu, anv@email.unc.edu}
\end{center}
\centerline{October, 2002}
\begin{abstract}
We present an elliptic version of Selberg's integral formula.
\end{abstract}
\section{Introduction}
The Selberg integral is the integral
\begin{equation*}
B_{p}(\alpha,\beta,\gm)=
\int_{\Delta_p}\prod_{j=1}^{p}t_{j}\,^{\alpha-1}(1-t_{j})^{\beta-1}
\prod_{0\leq j<k\leq 1}(t_{j}-t_{k})^{2\gm}dt_{1}\dots dt_{p},
\end{equation*}
where $\Delta_{p}=\{t\in\mathbb{R}^{p}\,|\,0\leq t_{p}\leq\dots\leq t_{1}\leq 1\}.$
The Selberg integral is a generalization of the beta function.  It can
be calculated explicitly,
\begin{equation*}
B_{p}(\alpha,\beta,\gm)=
{1\over p!}
\prod_{j=0}^{p-1}
\frac{\Gm(1+\gm+j\gm)
\Gm(\alpha+j\gm)
\Gm(\beta+j\gm)}
{\Gm(1+\gm)
\Gm(\alpha+\beta+(p+j-1)\gm)}.
\end{equation*}
The Selberg integral has many applications, see \cite{A1,A2,As,D,DF1,DF2,M,S}. In this paper, we
present elliptic versions of the Selberg integral.  
\section{Conformal Blocks on the Torus}
Let $\tau\in\mathbb{C}$ be such that $\Im\tau>0$.  Let $\ka$ and $p$ be
non-negative integers satisfying $\ka\geq 2p+2$.
The KZB-heat equation is the partial differential equation
\begin{equation}\label{KZB}
2\pi i\ka{\partial u\over\partial\tau}(\la,\tau)={\partial^2
  u\over\partial\la^{2}}(\la,\tau)+p(p+1)\rho\,'(\la,\tau)u(\la,\tau).
\end{equation}
Here, the prime denotes the derivative with respect to the first
argument, and $\rho$ is defined in terms of 
the first Jacobi theta function,
\begin{equation*}
\thi(\la,\tau)=2q^{\frac{1}{8}}
\sin{(\pi\la)}
\prod_{j=1}^{\infty}
(1-q^{j}e^{2\pi i\la})
(1-q^{j}e^{-2\pi i\la})
(1-q^{j})
,\quad
\rho(\la,\tau)=\frac{\thi'(\la,\tau)}{\thi(\la,\tau)}, 
\end{equation*}
where $q=e^{2\pi i\tau}$.  
Holomorphic solutions of the KZB-heat equation with the properties,
\begin{enumerate}
\renewcommand{\theenumi}{\roman{enumi}}\renewcommand{\labelenumi}{(\theenumi)}
\item $u(\la+2,\tau)=u(\la,\tau)$,
\item $u(\la+2\tau,\tau)=e^{-2\pi i\ka(\la+\tau)}\,u(\la,\tau)$, 
\item $u(-\la,\tau)=(-1)^{p+1}\,u(\la,\tau)$,
\item $u(\la,\tau)=\mathcal{O}((\la-m-n\tau)^{p+1})$ as $\la\to
  m+n\tau$ for any $m,n\in\mathbb{Z}$
\end{enumerate}
are called conformal blocks (or elliptic hypergeometric functions) associated with the family of elliptic curves
$\mathbb{C}/\mathbb{Z}+\tau\mathbb{Z}$ with the marked point $z=0$ and the
irreducible $sl_{2}$ representation of dimension $2p+1$.  It is known that the space of conformal blocks has dimension $\ka-2p-1$.   
\section{Integral Representations of Conformal Blocks}
Introduce special functions 
\begin{equation*}
\sigma_{\la}(t,\tau)=
\frac{\thi(\la-t,\tau)\thi'(0,\tau)}
{\thi(\la,\tau)\thi(t,\tau)},
\quad
E(t,\tau)=\frac{\thi(t,\tau)}{\thi'(0,\tau)}.
\end{equation*}
Consider the theta functions 
\begin{equation*}
\theta_{\ka,n}(\la,\tau)=
\sum_{j\in\Z}e^{2\pi i\ka(j+\frac{n}{2\ka})^{2}\tau+2\pi
  i\ka(j+\frac{n}{2\ka})\la},\quad n\in\Z/2\ka\Z.
\end{equation*}
They form a basis of the space of theta functions of level $\ka$. 
They satisfy the equations
\begin{gather*}
\theta_{\ka,n}(\la+1,\tau)=(-1)^{n}\theta_{\ka,n}\lt,\quad
\theta_{\ka,n}(\la+\tau,\tau)=e^{-\pi i\ka(\la+\frac{\tau}{2})}\theta_{\ka,n+\ka}\lt
\intertext{and have the modular properties}
\theta_{\ka,n}(\la,\tau+1)=e^{\pi i\frac{n^2}{2\ka}}\theta_{\ka,n}\lt,\quad
\theta_{\ka,n}\left({\la\over\tau},-{1\over\tau}\right)=\sqrt{-{i\tau\over2\ka}}e^{\pi
    i\ka{\la^{2}\over2\tau}}\,\sum_{m=0}^{2\ka-1}e^{-\pi i\frac{mn}{\ka}}\,\theta_{\ka,m}(\la,\tau),
\end{gather*}
where $|\text{arg}(-i\tau)|<\pi/2$.
Let $\tht{\ka}{n}$ denote the symmetrization of $\theta_{\ka,n}$ with
respect to $\la$,
$
\tht{\ka}{n}(\la,\tau)=\theta_{\ka,n}(\la,\tau)+\theta_{\ka,n}(-\la,\tau)$.

Define $u_{\ka,n}$ by 
\begin{equation*}
u_{\ka,n}(\la,\tau)=u_{p,\ka,n}(\la,\tau)=J_{p,\ka,n}(\la,\tau)+(-1)^{p+1}J_{p,\ka,n}(-\la,\tau),
\end{equation*}
where
\begin{multline*}
J_{p,\ka,n}(\la,\tau)
=
\int_{\Delta_{p}} 
\prod_{j=1}^p E(t_{j},\tau)^{-\frac{2p}{\ka}}
\prod_{1\leq j<k\leq
  p}E(t_{j}-t_{k},\tau)^{\frac{2}{\ka}}\,\times
\\
\prod_{j=1}^p\sigma_{\la}(t_{j},\tau)
\theta_{\ka,n}\left(\la+\frac{2}{\ka}\sum_{j=1}^p
t_{j},\tau\right)dt_{1}\dots dt_{p}\,\,.  
\end{multline*}
The branch of the logarithm is chosen in such a way that
arg $(E(t,\tau))\to 0$ as $t\to 0^{+}$, and the integral is understood
as the analytic continuation from the region where all of the exponents in the integrand have positive real parts. 
\begin{thm}\cite{FV1}\label{t0}
For all $n$, the integrals $u_{\ka,n}(\la,\tau)$ are solutions of the KZB-heat
equation having the properties (i)-(iv).    
\end{thm}
\begin{thm}\cite{FSV2}\label{t00}
We have
\begin{enumerate}
\renewcommand{\theenumi}{\alph{enumi}}\renewcommand{\labelenumi}{(\theenumi)}
\item 
$u_{\ka,n}=u_{\ka,n+2\ka}$ and
$u_{\ka,n}=-e^{2\pi ipn/\ka}u_{\ka,-n}$.
\item
The set
$
\{u_{\ka,n}\lt\,|\,n=p+1,\dots,\ka-p-1\}
$
is a basis for the space of conformal blocks.  The integrals $u_{\ka,n}$ are identically zero for all other values of $n$ in the interval from $0$ to $\ka$.
\end{enumerate}
\end{thm}

\section{Transformations acting on the space of conformal blocks}
Introduce four transformations $A$, $B$, $T$, and $S$ defined by
\begin{gather*}
Au(\la,\tau)=u(\la+1,\tau),\quad
Bu(\la,\tau)=e^{\pi i \ka(\la+\frac{\tau}{2})}u(\la+\tau,\tau),\\
Tu(\la,\tau)=u(\la,\tau+1),\quad
Su(\la,\tau)=e^{-\pi i \ka{\la^{2}\over2\tau}}\tau^{-\frac{1}{2}-{p(p+1)\over\ka}}\,u\left({\la\over\tau},-{1\over\tau}\right),
\end{gather*}
where we fix $\arg\tau\in(0,\pi)$.  

\begin{proposition}\label{pp1}
If $u(\la,\tau)$ is a solution of the KZB-heat equation, then
$Au(\la,\tau)$, $Bu\lt$, $Tu(\la,\tau)$, and $Su(\la,\tau)$
are solutions too.  Moreover, the transformations $A$, $B$, $T$ and $S$ preserve
the properties (i)-(iv).  
\end{proposition}
The proofs that $T$ and $S$ preserve the space of conformal blocks are given in \cite{EK}.  The proofs that $A$ and $B$ also preserve this space are straightforward and follow from the equations
\begin{gather*}
\thi(\la+1,\tau)=-\thi\lt,\quad
\thi(\la+\tau,\tau)= -e^{-\pi i(2\la+\tau)}\thi\lt.
\end{gather*}
\begin{lemma}
Restricted to the space of conformal blocks, the transformations $A$, $B$, $T$, and $S$ satisfy the relations 
\begin{gather*}
A^2=I,\quad B^2=I,\quad S^2=(-1)^p i e^{-\pi i\frac{p(p+1)}{\ka}}I,
\quad(ST)^{3}=(-1)^p i e^{-\pi i\frac{p(p+1)}{\ka}}I,\\
SAS^{-1}=B,\quad AB=(-1)^{\ka}BA,\quad TB=i^{\ka}BAT,
\end{gather*}
where $I$ denotes the identity transformation.\qed
\end{lemma}
\begin{lemma}\label{l1}
We have
\begin{gather*}
Au_{\ka,n}(\la,\tau)=(-1)^{n}u_{\ka,n}\lt,\quad
Bu_{\ka,n}\lt=-e^{2\pi i\frac{pn}{\ka}}u_{\ka,\ka-n}\lt.
\qed
\end{gather*}
\end{lemma}
Let $(t_{m,n})$ and $(s_{m,n})$ be
the matrices of the transformations $T$ and $S$, respectively, with
respect to the basis
$
 \{u_{\ka,n}(\la,\tau)\,\vert\, p+1\leq n\leq \ka -p-1\}
$, 
namely,
$
Tu_{\ka,n}=\sum_{m=p+1}^{\ka-p-1}t_{m,n}u_{\ka,m},\,
Su_{\ka,n}=\sum_{m=p+1}^{\ka-p-1}s_{m,n}u_{\ka,m}
$.
In Theorem \ref{fsvthm}, we give formulas for the matrices of $T$ and $S$ in terms of Macdonald polynomials of type $A_{1}$.

The Macdonald polynomials \cite{Ma} of type $A_{1}$ are $x$-even polynomials
in terms of $\epsilon^{mx}$, where $m\in\mathbb{Z}$.  They depend on two
parameters $k$ and $n$, where $k$ and $n$ are non-negative
integers.  They are
defined by the conditions:
\begin{enumerate}  
\item 
$P_{n}^{(k)}(x)=\epsilon^{nx}+\epsilon^{-nx}+$ lower order terms, except for $P_{0}^{(k)}(x)=1$,     
\item
$\langle P_{m}^{(k)},P_{n}^{(k)}\rangle=0$ for $m\neq n$, where
\begin{equation*}
\langle f,g\rangle= \,\frac{1}{2}\text{Const Term}\,\left(fg\prod_{j=0}^{k-1}\left(1-\epsilon^{2(j+x)}\right)\left(1-\epsilon^{2(j-x)}\right)\right).
\end{equation*}
\end{enumerate}
\begin{thm}\cite{FSV2}\label{fsvthm}
Let $\epsilon=e^{\pi i/\ka}$.
For $p+1\leq m,n\leq \ka-p-1$, we have
\begin{gather*}
t_{m,n}=\epsilon^{\frac{n^2}{2}}\delta_{mn},\\
s_{m,n}=
\frac{e^{-\frac{\pi i}{4}}}{\sqrt{2\ka}}\epsilon^{p(n-m)-\frac{p(p+1)}{2}}(\epsilon^{-m}-\epsilon^{m})\left(\prod_{j=1}^{p}(\epsilon^{-n+j}-\epsilon^{n-j})\right)P^{(p+1)}_{n-p-1}(m),
\end{gather*}
where $\delta_{mn}=1$ for $m=n$ and $0$ otherwise.
\end{thm}
\section{Integral identities}
To formulate our main result, we need functions $\eta$, $\phi_{1}$, $\phi_{2}$, and $\phi_3$.  
The Dedekind $\eta$-function is the function 
$
\eta(\tau)=
q^{1/24}
\prod_{j=1}^{\infty}
(1-q^{j})
$.
We have
$
\thi'(0,\tau)=2\pi \eta^{3}(\tau)
$.
Consider the functions \cite{W}
\begin{gather*}
\phi_{1}(\tau)=\frac{\eta(\tau)^{2}}{\eta(\frac{\tau}{2})\eta(2\tau)}=
q^{-\frac{1}{48}}\prod_{j=1}^{\infty}\left(1+q^{j-\frac{1}{2}}\right),\quad
\phi_{2}(\tau)=\frac{\eta(\frac{\tau}{2})}{\eta(\tau)}=
q^{-\frac{1}{48}}\prod_{j=1}^{\infty}\left(1-q^{j-\frac{1}{2}}\right),\\
\phi_{3}(\tau)=\sqrt{2}\,\frac{\eta(2\tau)}{\eta(\tau)}
=\sqrt{2}q^{\frac{1}{24}}\prod_{j=1}^{\infty}(1+q^{j}).
\end{gather*}
We have
\begin{gather*}
\phi_{1}\left(-\frac{1}{\tau}\right)=\phi_{1}(\tau),\quad
\phi_{2}\left(-\frac{1}{\tau}\right)=\phi_{3}(\tau),\quad
\phi_{3}\left(-\frac{1}{\tau}\right)=\phi_{2}(\tau),\\
\phi_{1}(\tau+1)=e^{-\frac{\pi i}{24}}\phi_{2}(\tau),\quad
\phi_{2}(\tau+1)=e^{-\frac{\pi i}{24}}\phi_{1}(\tau),\quad
\phi_{3}(\tau+1)=e^{\frac{\pi i}{12}}\phi_{3}(\tau).
\end{gather*}
Let
\begin{equation*}
c_{\ka,n}=c_{p,\ka,n}=(2\pi)^{\frac{p(p+1)}{\ka}}e^{-\pi
  i\frac{p(3p-1)}{2\ka}}e^{\pi i\frac{p+1}{2}}
B_{p}\left(\frac{n+1}{\ka},-\frac{2p}{\ka},\frac{1}{\ka}\right)
\prod_{j=1}^{p}\left(1-e^{2\pi i\frac{n+j}{\ka}}\right).
\end{equation*}
Here, $B_{p}(\alpha,\beta,\gamma)$ is the Selberg integral.
\begin{thm}\label{t1}
We have ten series of identities,
\begin{gather}\label{esi1} 
 u_{2p+2,p+1}(\la,\tau)=c_{2p+2,p+1}
\thi(\la,\tau)^{p+1},\\\label{esi2}
u_{2p+3,p+1}\lt=c_{2p+3,p+1}
\eta(\tau)^{-\frac{3(p+1)}{2p+3}}\thi^{p+1}\lt\theta_{1,0}\lt,\\\label{esi3}
u_{2p+3,p+2}\lt=c_{2p+3,p+2}
\eta(\tau)^{-\frac{3(p+1)}{2p+3}}\thi^{p+1}\lt\theta_{1,1}\lt,\\\label{esi4}
 u_{2p+4,p+2}\lt=2^{-\frac{2(p+1)}{2p+4}}c_{2p+4,p+2}
\left(\phi_{3}(\tau)\eta(\tau)^{-1}\right)^{\frac{4(p+1)}{2p+4}}
\thi^{p+1}\lt\tht{2}{1}\lt,
\end{gather} 
\vspace{-.3 in}
\begin{multline}\label{esi5}
u_{2p+4,p+1}\lt+(-1)^{p+1}e^{2\pi i\frac{p(p+1)}{2p+4}}u_{2p+4,p+3}\lt=\\
c_{2p+4,p+1}
\left(\phi_{2}(\tau)\eta(\tau)^{-1}\right)^{\frac{4(p+1)}{2p+4}}
\thi^{p+1}\lt
(\theta_{2,0}\lt - \theta_{2,2}\lt)
,\end{multline}
\vspace{-.28 in}
\begin{multline}\label{esi6}
u_{2p+4,p+1}\lt+(-1)^{p}e^{2\pi i\frac{p(p+1)}{2p+4}}u_{2p+4,p+3}\lt=\\
c_{2p+4,p+1}
\left(\phi_{1}(\tau)\eta(\tau)^{-1}\right)^{\frac{4(p+1)}{2p+4}}
\thi^{p+1}\lt
(\theta_{2,0}\lt + \theta_{2,2}\lt),\end{multline}
\vspace{-.265 in}
\begin{multline}
\label{esi7}
u_{2p+6,p+1}\lt+(-1)^{p+1}e^{2\pi i\frac{p(p+1)}{2p+6}}u_{2p+6,p+5}\lt=\\
2^{\frac{3(p+1)}{2p+6}}c_{2p+6,p+1}
\left(\phi_{3}(\tau)\eta(\tau)\right)^{-\frac{6(p+1)}{2p+6}}
\thi^{p+1}\lt
(\theta_{4,0}\lt-\theta_{4,4}\lt)
,
\end{multline}
\vspace{-.265 in}
\begin{multline}
\label{esi8}
u_{2p+6,p+2}\lt+(-1)^{p}e^{2\pi i\frac{p(p+2)}{2p+6}}u_{2p+6,p+4}\lt=\\
c_{2p+6,p+2} 
\left(\phi_{2}(\tau)\eta(\tau)\right)^{-\frac{6(p+1)}{2p+6}}
\thi^{p+1}\lt
(\theta_{4,1}^{s}\lt+\theta_{4,3}^{s}\lt)
,
\end{multline}
\vspace{-.265 in}
\begin{multline}
\label{esi9}
u_{2p+6,p+2}\lt+(-1)^{p+1}e^{2\pi i\frac{p(p+2)}{2p+6}}u_{2p+6,p+4}\lt=\\
c_{2p+6,p+2}
\left(\phi_{1}(\tau)\eta(\tau)\right)^{-\frac{6(p+1)}{2p+6}}
\thi^{p+1}\lt
(\theta_{4,1}^{s}\lt-\theta_{4,3}^{s}\lt),
\end{multline}
\vspace{-.265 in}
\begin{multline}
\label{esi10}
u_{2p+8,p+2}\lt+(-1)^{p+1}e^{2\pi i\frac{p(p+2)}{2p+8}}u_{2p+8,p+6}\lt=\\
c_{2p+8,p+2}
\eta(\tau)^
{\frac{-8(p+1)}{2p+8}}
\vartheta_{1}^{p+1}\lt
(\theta_{6,1}^{s}\lt-\theta_{6,5}^{s}\lt)
.
\end{multline}
\end{thm}
The integrals in Theorem \ref{t1} are
appropriately called the elliptic Selberg integrals.  The identity \Ref{esi1} appears in \cite{FSV1} and in \cite{FV1} for $p=1$.
\section{Differential equations} 
In Lemmas \ref{del1}--\ref{del3}, let $'$ denote the derivative with respect to $\la$, let $\dot{}$ denote the derivative with respect to $\tau$, and let
$
v(\la,\tau)=\thi^{p+1}\lt\sum_{j=0}^{\ka-2p-2}c_{j}(\tau)\theta^{s}_{\ka-2p-2,j}\lt$.
\begin{lemma}\label{del1}
The function $v(\la,\tau)$
is a solution of the KZB-heat equation if and only if the differential equation
\begin{multline}\label{diffeq} 
\frac{\ka}{p+1}\sum_{j=0}^{\ka-2p-2}\left(\frac{d}{d\tau}c_{j}\right)\theta^{s}_{\ka-2p-2,j}=
(2p+2-\ka)\frac{\dot{\thi}}{\thi}\sum_{j=0}^{\ka-2p-2}c_{j}\theta^{s}_{\ka-2p-2,j}\\
-2\sum_{j=0}^{\ka-2p-2}c_{j}\dot{\theta}^{s}_{\ka-2p-2,j}+
\frac{1}{\pi i}\frac{\thi'}{\thi}\sum_{j=0}^{\ka-2p-2}c_{j}(\theta^{s}_{\ka-2p-2,j})'
\end{multline}
holds.  
\end{lemma}
The proof of Lemma \ref{del1} 
uses the identities
\begin{gather}\label{e1}
2\pi
i(2p+2)\dot{(\vartheta_{1}^{p+1})}(\la,\tau)=(\thi^{p+1})''(\la,\tau)
+p(p+1)\rho'(\la,\tau)\thi^{p+1}(\la,\tau),\\\label{e2}
2\pi i\ka\dot{\theta}^{s}_{\ka,m}(\la,\tau)=(\theta^{s}_{\ka,m})''(\la,\tau).
\end{gather}
\begin{proof}[Proof of Lemma \ref{del1}]
Applying the differential operator
$2\pi i\ka\partial/\partial\tau$ to $v(\la,\tau)$ gives
\begin{multline*}
2\pi
i\ka\biggl[
(p+1)\dot{\thi}\thi^{p}\sum_{j=0}^{\ka-2p-2}c_{j}\theta^{s}_{\ka-2p-2,j}
+\thi^{p+1}\sum_{j=0}^{\ka-2p-2}
\left(\frac{d}{d\tau}c_{j}\right)\theta^{s}_{\ka-2p-2,j}\\
+\thi^{p+1}\sum_{j=0}^{\ka-2p-2}c_{j}\dot{\theta}^{s}_{\ka-2p-2,j}\biggl].
\end{multline*}
Applying the differential operator 
$\partial^{2}/\partial\la^{2}+p(p+1)\rho'(\la,\tau)$ to
$v(\la,\tau)$ gives
\begin{multline*}
(\thi^{p+1})''\sum_{j=0}^{\ka-2p-2}c_{j}\theta^{s}_{\ka-2p-2,j}+2(p+1)\thi'\thi^{p}\sum_{j=0}^{\ka-2p-2}c_{j}(\theta^{s}_{\ka-2p-2,j})'\\
+\thi^{p+1}\sum_{j=0}^{\ka-2p-2}c_{j}(\theta^{s}_{\ka-2p-2,j})''+p(p+1)\rho'\thi^{p+1}\sum_{j=0}^{\ka-2p-2}c_{j}\theta^{s}_{\ka-2p-2,j}.
\end{multline*}
Applying \Ref{e1} and \Ref{e2}, we obtain the result.
\end{proof}
\begin{lemma}\label{del}  If  $v(\lambda,\tau)$
is a solution of the KZB-heat equation, then the functions $c_{j}(\tau)$ satisfy the differential equation 
\begin{multline*}
\frac{\ka}{p+1}\sum_{j=0}^{\ka-2p-2}\left(\frac{d}{d\tau}c_{j}(\tau)\right)\theta_{\ka-2p-2,j}(0,\tau)=\\
(2p+2-\ka)\left(\frac{d}{d\tau}\ln\vartheta_{1}'(0,\tau)\right)
\sum_{j=0}^{\ka-2p-2}c_{j}(\tau)\theta_{\ka-2p-2,j}(0,\tau)\\
+
2(\ka-2p-3)
\sum_{j=0}^{\ka-2p-2}c_{j}(\tau)
\frac{d}{d\tau}\theta_{\ka-2p-2,j}(0,\tau).
\end{multline*}
\end{lemma}
\begin{proof}
For any fixed $\la$, equation \Ref{diffeq} gives a differential equation for the functions $c_{j}(\tau)$.  We take the limit of that equation as $\la\to 0$.  In the
ratio $\dot{\thi}/\thi$, both the numerator and
denominator tend to zero, so the limit of this term as $\la\to 0$ is
equal to the limit of the ratio of the derivatives of the numerator and the denominator.   
The limit of the ratio $(\sum_{j=0}^{\ka-2p-2}c_{j}(\theta^{s}_{\ka-2p-2,j})')/\thi$ is calculated in
the same way, since each $\theta^{s}_{\ka-2p-2,j}$ is a symmetric function and therefore $(\theta^{s}_{\ka-2p-2,j})'(0,\tau)=0$. 
Then the result follows from \Ref{e2} .
\end{proof}
\begin{lemma}\label{del3}
If $v(\la,\tau)$
is a solution of the KZB-heat equation, then the functions $c_{j}(\tau)$ satisfy the differential equation 
\begin{multline*}
\frac{\ka}{p+1}\sum_{j=0}^{\ka-2p-2}\left(\frac{d}{d\tau}c_{j}(\tau)\right)\theta_{\ka-2p-2,\ka-2p-2-j}(0,\tau)=\\
(2p+2-\ka)\left(\frac{d}{d\tau}\ln\thi'(0,\tau)\right)
\sum_{j=0}^{\ka-2p-2}c_{j}(\tau)\theta_{\ka-2p-2,\ka-2p-2-j}(0,\tau)\\
+
2(\ka-2p-3)
\sum_{j=0}^{\ka-2p-2}c_{j}(\tau)
\frac{d}{d\tau}\theta_{\ka-2p-2,\ka-2p-2-j}(0,\tau).
\end{multline*}
\end{lemma}
\begin{proof}
We take the limit of \Ref{diffeq} as $\la\to\tau$.  This limit is calculated in terms of the limit $\la\to0$ using the formulas
\begin{gather*}
\frac{\partial}{\partial z}\theta_1(z, \tau)\vert_{z=\la+\tau} = 
e^{- 2\pi i \lambda - \pi i \tau}\left(2\pi i\theta_1(\lambda, \tau)
- \theta_1'(\lambda, \tau)\right),\\
\frac{\partial}{\partial z} \theta_1(\lambda + \tau, z)\vert_{z=\tau}  = 
e^{- 2\pi i \lambda - \pi i \tau}\left(-\pi i\theta_1(\lambda, \tau)
+ \theta_1'(\lambda, \tau)
- \dot \theta_1(\lambda, \tau)\right),\\
\frac{\partial}{\partial z}\theta_{\kappa,m}^{s}(z, \tau)\vert_{z=\la+\tau} = 
e^{- \pi i \kappa\lambda - \pi i\kappa{ \tau\over2}}\left(-\pi i\ka\theta_{\kappa,\kappa-m}^{s}(\lambda, \tau)+
(\theta_{\kappa,\kappa-m}^{s})'(\lambda, \tau)\right),\\
\frac{\partial}{\partial z}\theta_{\kappa,m}^{s}(\lambda + \tau, z)\vert_{z=\tau}  = 
 e^{-\pi i\kappa \lambda - \pi i\kappa{ \tau\over2}}
\left(\pi i{\kappa\over2}\theta_{\kappa,\kappa-m}^{s}(\lambda, \tau)
-(\theta_{\kappa,\kappa-m}^{s})'(\lambda, \tau)
+
\dot \theta_{\kappa,\kappa-m}^{s}(\lambda, \tau)\right).
\end{gather*}
It is straightforward to calculate the limit of the left hand side.  Using the above formulas, the limit as $\la\to \tau$ of the right hand side is equal to the limit as $\la\to 0$ of the expression
\begin{multline*}
e^{\pi i(2p+2-\ka)(\la+\frac{\tau}{2})}
\biggl((2p+2-\ka)\frac{\dot{\thi}}{\thi}
\sum_{j=0}^{\ka-2p-2}c_{j}\theta^{s}_{\ka-2p-2,\ka-2p-2-j}\\
-2\sum_{j=0}^{\ka-2p-2}c_{j}\dot{\theta}^{s}_{\ka-2p-2,\ka-2p-2-j}+
\frac{1}{\pi i}\frac{\thi'}{\thi}
\sum_{j=0}^{\ka-2p-2}c_{j}(\theta^{s}_{\ka-2p-2,\ka-2p-2-j})'
\biggl).
\end{multline*}
This limit is calculated using L'H\^{o}pital's rule.
\end{proof} 
\section{Identities for theta functions}
In the next section, we give the proofs of the integral identities in 
Theorem \ref{t1}. 
We will use the following results.
\begin{lemma}\label{p1} We have
$
\tht{2}{1}(\la)=\thi(\la+1/2)
$.
\end{lemma} 
Lemma \ref{p1} is proved by comparing the Fourier series expansions of the functions.
\begin{corollary}\label{c1}
We have  
$2\theta_{2,1}(0)=\eta(\tau)\phi_{3}(\tau)^{2}.$
\end{corollary}
\begin{lemma}\label{p3}
Let 
\begin{gather*}
f_{1}(\tau)=\frac{\theta_{4,1}(0)-\theta_{4,3}(0)}{\eta(\tau)},\quad
f_{2}(\tau)=\frac{\theta_{4,1}(0)+\theta_{4,3}(0)}{\eta(\tau)},\quad
f_{3}(\tau)=\frac{\theta_{4,0}(0)-\theta_{4,4}(0)}{\sqrt{2}\eta(\tau)}.
\end{gather*}
Then
$
f_{1}(\tau)=\phi_{1}(\tau)^{-1},\, 
f_{2}(\tau)=\phi_{2}(\tau)^{-1},\,
f_{3}(\tau)=\phi_{3}(\tau)^{-1}.
$
\end{lemma}
The proof of Lemma \ref{p3} is based on the following result.
\begin{lemma}\label{wa}\cite{W}
Suppose $g_{1}(\tau)$, $g_{2}(\tau)$, and $g_{3}(\tau)$ are holomorphic functions on the upper half plane $\mathbb{C}_{+}$ satisfying the following conditions.
\begin{enumerate}
\item[P1.] The functions $g_{1}$, $g_{2}$, and $g_{3}$  can be written in the forms
\begin{gather*}
g_{1}(\tau)=q^{-\frac{a}{48}}\sum_{j=0}^{\infty}a_{j}q^{\frac{j}{2}},\quad
g_{2}(\tau)=q^{-\frac{a}{48}}\sum_{j=0}^{\infty}(-1)^{j}a_{j}q^{\frac{j}{2}},\quad
g_{3}(\tau)=q^{\frac{a}{24}}\sum_{j=0}^{\infty}b_{j}q^{j},
\end{gather*}
where $a$ is an integer, $a_{j},b_{j}\in\mathbb{C}$, and $a_{0}=1$.
\item[P2.] The functions $g_{1}$, $g_{2}$, and $g_{3}$ have the modular properties
\begin{gather*}
g_{1}\left(-\frac{1}{\tau}\right)=g_{1}(\tau),\quad
g_{2}\left(-\frac{1}{\tau}\right)=g_{3}(\tau),\quad 
g_{3}\left(-\frac{1}{\tau}\right)=g_{2}(\tau).
\end{gather*}
Then $g_{i}(\tau)=\phi_{i}(\tau)^{a}$, for $i=1,2,3$.
\end{enumerate}
\end{lemma}
\begin{proof}[Proof of Proposition \ref{p3}]
We have
\begin{gather*}
\theta_{4,1}(0)=
\sum_{j\in\mathbb{Z}}q^{4(j+\frac{1}{8})^{2}}
=q^{\frac{1}{16}}\sum_{j\in\mathbb{Z}}q^{\frac{1}{2}(8j^{2}+2j)},
\quad
\theta_{4,3}(0)=
\sum_{j\in\mathbb{Z}}q^{4(j+\frac{3}{8})^{2}}=
q^{\frac{1}{16}}\sum_{j\in\mathbb{Z}}q^{\frac{1}{2}(8j^{2}+6j+1)},
\\
\theta_{4,0}(0)=\sum_{j\in\mathbb{Z}}q^{4j^{2}}=
\sum_{j\in\mathbb{Z}}q^{(2j)^{2}},
\quad
\theta_{4,4}(0)=\sum_{j\in\mathbb{Z}}q^{4(j+\frac{1}{2})^{2}}=
\sum_{j\in\mathbb{Z}}q^{(2j+1)^{2}}.
\end{gather*}
Hence, the functions 
\begin{gather*}
\tilde{f}_{1}(\tau)=q^{-\frac{1}{24}}(\theta_{4,1}(0)-\theta_{4,3}(0)),\quad
\tilde{f}_{2}(\tau)=q^{-\frac{1}{24}}(\theta_{4,1}(0)+\theta_{4,3}(0)),\\ 
\tilde{f}_{3}(\tau)=2^{-\frac{1}{2}}q^{-\frac{1}{24}}(\theta_{4,0}(0)-\theta_{4,4}(0))
\end{gather*}
are holomorphic functions on $\mathbb{C}_{+}$ with the property P1 for $a=-1$ and $g_{i}=\tilde{f}_{i}$.
Let $y(q)=\sum_{j=1}^{\infty}c_{j}q^{j}$ be defined by the condition
$
1+y(q)=q^{-1/24}\eta(\tau).
$
Then for $i=1,2,3$, 
\begin{equation*}
f_{i}(\tau)=\frac{\tilde{f}_{i}(\tau)}{1+y(q)}
=\tilde{f}_{i}(\tau)(1-y(q)+y(q)^{2}+\dots)
\end{equation*}
are holomorphic functions on $\mathbb{C}_{+}$ with the property P1 for $a=-1$ and $g_{i}={f}_{i}$.
One checks that $f_{1}$, $f_{2}$ and $f_{3}$ have the property P2 using the modular properties of $\theta_{4,n}(0)$ and $\eta(\tau)$.
\end{proof}
\begin{lemma}\label{p4}
We have
$
\theta_{6,1}(0)-\theta_{6,5}(0)=\eta(\tau).
$
\end{lemma}
Lemma \ref{p4} is proved by comparing the infinite series expansions of the functions. 
\section{The proof of Theorem \ref{t1}}
\subsection{Proof of \Ref{esi1}}
For $\ka=2p+2$,
the space of conformal blocks is one-dimensional.
The right hand side of \Ref{esi1} is a solution
of \Ref{KZB} with the properties (i)-(iv) \cite{FV1}.  According to
Theorem \ref{t0}, the left hand side also has these properties.  Thus the
two functions are proportional.  The coefficient of proportionality is
calculated by comparing the leading terms of $\thi^{p+1}$ and
$u_{\ka,p+1}$ in the limit as $\tau\to i\infty$.  The leading term of
$\thi^{p+1}$ is
$
(-i)^{p+1}q^{(p+1)/8}(e^{\pi i\la}-e^{-\pi i\la})^{p+1}$.
Let $dt=dt_{1}\dots dt_{p}$.  The leading term of $u_{\ka,p+1}$ is
\begin{multline*}
\int_{\Delta_{p}}
\prod_{j=1}^{p}
\left(
\frac{e^{\pi it_{j}}-e^{-\pi it_{j}}}{2\pi e^{\frac{\pi i}{2}}}
\right)
^{-\frac{2p}{2p+2}-1}
\prod_{1\leq j<k\leq p}
\left(
\frac{e^{\pi i(t_{j}-t_{k})}-e^{-\pi i(t_{j}-t_{k})}}{2\pi
  e^{\frac{\pi i}{2}}}
\right)
^{\frac{2}{2p+2}}\\
\biggl(
\left(
\prod_{j=1}^{p}
\frac{e^{\pi i(\la-t_{j})}-e^{-\pi i(\la-t_{j})}}{e^{\pi i\la}-e^{-\pi
    i\la}}
\right)
q^{\frac{(p+1)^{2}}{4(2p+2)}}
e^{\pi i(p+1)
(\la+\frac{2}{2p+2}\sum_{j=1}^{p}t_{j})
}\\
+(-1)^{p+1}
\left(\prod_{j=1}^{p}\frac{e^{\pi i(\la+t_{j})}-e^{-\pi i(\la+t_{j})}}
{e^{\pi i\la}-e^{-\pi i\la}}\right)
q^{\frac{(p+1)^{2}}{4(2p+2)}}
e^{\pi i(p+1)(-\la+\frac{2}{2p+2}\sum_{j=1}^{p}t_{j})}
\biggl)
dt.
\end{multline*}
The above expression is equal to
\begin{gather*}
(2\pi e^{\frac{\pi i}{2}})^{\frac{p(p+1)}{2p+2}+p}e^{-\pi i(\frac{2p^{2}}{2p+2}+p)}
   q^{\frac{p+1}{8}}
(e^{\pi i\la}-e^{-\pi i\la})^{-p}
\sum_{l=0}^{p}I_{l}\,(e^{\pi i(2l+1)\la}-e^{-\pi i(2l+1)\la}),
\intertext{where}
I_{l}=
\int_{\Delta_{p}}
f_{l}(t_{1},\dots,t_{p})
\prod_{j=1}^{p}
e^{2\pi it_{j}(\frac{p+2}{2p+2}+\half)}
(1-e^{2\pi it_{j}})^{-\frac{2p}{2p+2}-1}
\prod_{1\leq j<k\leq p}(e^{2\pi it_{j}}-e^{2\pi
    it_{k}})^{\frac{2}{2p+2}}
dt,
\end{gather*}
for some function $f_{l}$, symmetric in the variables $t_{1},\dots,t_{p}$. 
Comparing the coefficients of $e^{\pi i (p+1)\la}$ in the leading terms, we find that 
$
u_{\ka,p+1}=i^{p+1}(2\pi e^{\frac{\pi i}{2}})^{\frac{p(p+1)}{2p+2}+p}e^{-\pi i(\frac{2p^{2}}{2p+2}+p)}I_{p}\thi^{p+1}$.
To complete the proof, it remains to compute $I_{p}$.
It is not difficult to show that
$
f_{p}(t_{1},\dots,t_{p})=\prod_{j=1}^{p}e^{-\pi it_{j}}$.
Let $x_{j}=e^{2\pi it_{j}}$.  Let $\tilde{\Delta}_{p}$ be the image of
$\Delta_{p}$ under the map $t_{j}\mapsto x_{j}$.  We have
\begin{equation*}
I_{p}=(2\pi i)^{-p}\int_{\tilde{\Delta}_{p}}
\prod_{j=1}^{p}x_{j}^{\frac{p+2}{2p+2}-1}(1-x_{j})^{-\frac{2p}{2p+2}-1}
\prod_{1\leq j<k\leq p}(x_{j}-x_{k})^\frac{2}{2p+2}dx.
\end{equation*}
Applying the Stokes theorem, we deform the contour
$\tilde{\Delta}_{p}$ to get
\begin{equation*}
I_{p}=(2\pi i)^{-p}\prod_{j=1}^{p}(e^{2\pi i\frac{j+p+1}{2p+2}}-1)
\int_{\Delta_{p}}
\prod_{j=1}^{p}x_{j}^{\frac{p+2}{2p+2}-1}(1-x_{j})^{-\frac{2p}{2p+2}-1}
\prod_{1\leq j<k\leq p}(x_{j}-x_{k})^\frac{2}{2p+2}dx.
\end{equation*}
Observe that
\begin{equation*}
\int_{\Delta_{p}}
\prod_{j=1}^{p}x_{j}^{\frac{p+2}{2p+2}-1}(1-x_{j})^{-\frac{2p}{2p+2}-1}
\prod_{1\leq j<k\leq p}(x_{j}-x_{k})^\frac{2}{2p+2}dx
\end{equation*}
is the Selberg integral $B_p((p+2)/(2p+2),-2p/(2p+2), 1/(2p+2))$.
This completes the proof.

\subsection{Proof of \Ref{esi2} and \Ref{esi3}}  
Let $\kappa=2p+3$.  Then any solution of \Ref{KZB} with the properties (i)-(iv) has the form $v\lt=\thi^{p+1}\lt(c_0(\tau)\theta_{1,0}\lt+c_1(\tau)\theta_{1,1}\lt)$.  Let $A$ be the transformation introduced in section $4$.  By Proposition \ref{pp1}, 
$
Av\lt=(-1)^{p+1}\thi^{p+1}\lt(c_0(\tau)\theta_{1,0}\lt-c_1(\tau)\theta_{1,1}\lt)
$ is also a solution.
Hence, for $j=0$ or $1$, the function  
$
v_{j}\lt=c_{j}(\tau)\thi^{p+1}\lt\theta_{1,j}\lt
$
gives a solution too.  Moreover, $Av_j=(-1)^{p+1+j}v_j$.
By Theorem \ref{t00}, the integrals $u_{\ka,p+1}$ and $u_{\ka,p+2}$
span the space of conformal blocks.  By Lemma \ref{l1},
$
Au_{\ka,p+1}=(-1)^{p+1}u_{\ka,p+1}$
and
$
Au_{\ka,p+2}=(-1)^{p}u_{\ka,p+2}$.
So for $j=0$ or $1$, the integral $u_{\ka,p+1+j}$ is proportional to
$v_{j}$.  
By Lemma \ref{del}, $c_{j}(\tau)$ must satisfy the differential equation
$$
\frac{\ka}{p+1}\left(\frac{d}{d\tau}c_{j}(\tau)\right)\theta_{1,j}(0,\tau)=-
\left(\frac{d}{d\tau}\ln\thi'(0,\tau)\right)c_{j}(\tau)\theta_{1,j}(0,\tau).$$
The function 
$
c_{j}(\tau)=\left((2\pi)^{-1}\thi'(0,\tau)\right)^{-(p+1)/\ka}=\eta(\tau)^{-3(p+1)/\ka}
$ 
is a solution of this equation for $j=0$ and $j=1$.
The coefficients of proportionality are computed in the limit
as $\tau\to i\infty$, cf. the proof of \Ref{esi1}.  This completes the proof.
\subsection{Proof of \Ref{esi4}, \Ref{esi5}, and \Ref{esi6}}
Let $\ka=2p+4$.  Let $v_{j}\lt=\thi^{p+1}\lt\tht{2}{j}\lt,\, 0\leq j\leq2$.
Any solution of \Ref{KZB} with the properties (i)-(iv) has the form $v\lt=\linebreak\sum_{j=0}^{2}c_j(\tau)v_j\lt$.
Let $A$ be the transformation in section $4$.  By Proposition \ref{pp1},
$Av=\sum_{j=0}^{2}(-1)^{p+j+1}c_jv_j$ is also a solution.
Hence the function $c_1v_1$ gives a solution too.  Moreover, it is an eigenvector of $A$ with eigenvalue $(-1)^p$.     
By Theorem \ref{t00}, the space of conformal blocks is three-dimensional 
with spanning set $\{u_{\ka,n}\,|\, p+1\leq n \leq
p+3\}$.  According to Lemma \ref{l1}, the eigenspace of $A$ corresponding to the eigenvalue $(-1)^{p}$ is one-dimensional and is spanned by $u_{\ka,p+2}$.
It follows that $u_{\ka,p+2}$ is proportional to $c_1 v_1$.  By Lemma \ref{del}, $c_1(\tau)$ must satisfy the differential equation
$$
\frac{\ka}{p+1}\left(\frac{d}{d\tau}c_{1}(\tau)\right)\theta_{2,1}(0,\tau)=-2\left(\frac{d}{d\tau}\ln\thi'(0,\tau)\right)c_{1}(\tau)\theta_{2,1}(0,\tau)+2c_{1}(\tau)\frac{d}{d\tau}\theta_{2,1}(0,\tau).$$
The function $c_{1}(\tau)=\left(4\pi\theta_{2,1}(0,\tau)\thi'(0,\tau)^{-1}\right)^{2(p+1)/\ka}$ is a solution of the above equation. 
By Corollary \ref{c1},
$
c_{1}(\tau)=\left(\phi_{3}(\tau)\eta(\tau)^{-1}\right)^{4(p+1)/\ka}$.
The coefficient of proportionality is computed in the limit $\tau\to i\infty$, cf. the proof of \Ref{esi1}.
This proves \Ref{esi4}.  
To prove \Ref{esi5}, we apply the transformation $S$ to both sides of \Ref{esi4}.  To prove \Ref{esi6}, we apply the transformation $T$ to both sides of \Ref{esi5}. 
This completes the proof. 

\subsection{Proof of \Ref{esi7}, \Ref{esi8}, and \Ref{esi9}}
Let $\ka=2p+6$.  Let $v_{j}\lt=\thi^{p+1}\lt\tht{4}{j}\lt,\, 0\leq j\leq4$.  Any solution of \Ref{KZB} with the properties (i)-(iv) has the form $v\lt=\newline\sum_{j=0}^{4}c_j(\tau)v_j\lt$.
Let $A$ and $B$ be the transformations in section $4$.  By Proposition \ref{pp1},
$Av=\sum_{j=0}^{4}(-1)^{p+j+1}c_jv_j$
is also a solution.
Hence the function $c_0v_0+c_2v_2+c_4 v_4$ gives a solution too.  Moreover,  
$B(c_0v_0+c_2v_2+c_4 v_4)=(-1)^{p+1}(c_4 v_0+c_2 v_2+c_0v_4)$ is also a solution.
So there exists a solution of the form $c(\tau)(v_0-v_4)$.  It is an eigenvector of $A$ with eigenvalue $(-1)^{p+1}$ and an eigenvector of B with eigenvalue $(-1)^p$. 
We show that the subspace of conformal blocks with this property is one-dimensional.
By Theorem \ref{t00}, the space of conformal blocks is five-dimensional 
with spanning set $\{u_{\ka,n}\,|\, p+1\leq n \leq
p+5\}$.   
By Lemma \ref{l1}, the eigenspace of $A$ corresponding to the eigenvalue $(-1)^{p+1}$ is three-dimensional and is spanned by $u_{\ka,p+1}$, $u_{\ka,p+3}$, and $u_{\ka,p+5}$.  By Lemma \ref{l1}, the transformation $B$ preserves the subspace $\langle u_{\ka,p+1},u_{\ka,p+3},u_{\ka,p+5}\rangle$.  The matrix of $B$ restricted to this subspace is 
$$
\begin{pmatrix} 0 & 0&-e^{2\pi i\frac{p(p+5)}{\ka}}\\
0&(-1)^{p+1}&0\\
-e^{2\pi i\frac{p(p+1)}{\ka}} & 0&0\end{pmatrix}.
$$
Thus the restriction of $B$ to  $\langle u_{\ka,p+1},u_{\ka,p+3},u_{\ka,p+5}\rangle$
has eigenvalues $(-1)^{p}$ and $(-1)^{p+1}$ of multiplicities $1$ and $2$, respectively.  The eigenspace corresponding to the eigenvalue $(-1)^{p}$ is spanned by the vector
$
u_{\ka,p+1}+(-1)^{p+1}e^{2\pi ip(p+1)/\ka}u_{\ka,p+5}.
$
It follows that this vector is proportional to 
$c(\tau)(v_{0}-v_{4})$. 
 By Lemma \ref{del}, $c(\tau)$ must satisfy the differential equation
\begin{multline*}
\frac{\ka}{p+1}\left(\frac{d}{d\tau}c(\tau)\right)(\theta_{4,0}(0,\tau)-\theta_{4,4}(0,\tau))=\\-4\left(\frac{d}{d\tau}\ln\thi'(0,\tau)\right)c(\tau)(\theta_{4,0}(0,\tau)-\theta_{4,4}(0,\tau))+6c(\tau)\frac{d}{d\tau}(\theta_{4,0}(0,\tau)-\theta_{4,4}(0,\tau)). 
\end{multline*}
The function $c(\tau)=\left((2\pi)^{2}(\theta_{4,0}(0,\tau)-\theta_{4,4}(0,\tau))^{3}\thi'(0,\tau)^{-2}\right)^{2(p+1)/\ka}$ is a solution of the above equation.
By Lemma \ref{p3}, we have
$c(\tau)=2^{3(p+1)/\ka}(\phi_{3}(\tau)\eta(\tau))^{-6(p+1)/\ka}$.
The coefficient of proportionality is computed as in the proof of \Ref{esi1}. 
This proves \Ref{esi7}.  To prove \Ref{esi8} and \Ref{esi9}, apply the transformations $S$ and $TS$, respectively, to both sides of \Ref{esi7}.  This completes the proof.
\subsection{Proof of \Ref{esi10}}
Let $\ka=2p+8$.  Let $v_{j}\lt=\thi^{p+1}\lt\tht{6}{j}\lt,\, 0\leq j\leq6.$  Any solution of \Ref{KZB} with the properties (i)-(iv) has the form $v\lt=\sum_{j=0}^{6}c_j(\tau)v_j\lt$.
Let $A$ and $B$ be as in section $4$.  By Proposition \ref{pp1}, $Av=\sum_{j=0}^{6}(-1)^{p+j+1}c_jv_j$ is also a solution.
Hence the function $c_1v_1+c_3 v_3+c_5 v_5$ gives a solution too.  The function
$
B(c_1v_1+c_3 v_3+c_5 v_5)=(-1)^{p+1}(c_5v_1+c_3v_3+c_1 v_5)
$
also gives a solution.  So there is a solution of the form $c(\tau)(v_1-v_5)$ which is an eigenvector of $A$ and $B$ with eigenvalue $(-1)^p$ under both transformations.  We show that the subspace of conformal blocks with this property is one-dimensional.  
By Theorem \ref{t00}, the space of conformal blocks is seven-dimensional 
with spanning set $\{u_{\ka,n}\,|\, p+1\leq n \leq
p+7\}$. 
By Lemma \ref{l1}, the three-dimensional eigenspace of $A$ corresponding to the eigenvalue $(-1)^{p}$ is spanned by $u_{\ka,p+2}$, $u_{\ka,p+4}$, and $u_{\ka,p+6}$.
The matrix of $B$ restricted to $\langle u_{\ka,p+2}, u_{\ka,p+4}, u_{\ka,p+6}\rangle$ is  
$$
\begin{pmatrix} 0&0 & -e^{2\pi i\frac{p(p+6)}{\ka}}\\
0&(-1)^{p+1}&0\\
-e^{2\pi i\frac{p(p+2)}{\ka}} &0& 0\end{pmatrix}.
$$
Thus, $B$ has eigenvalues $(-1)^{p}$ and $(-1)^{p+1}$ of multiplicities $1$ and $2$, respectively.  The eigenspace corresponding to the eigenvalue $(-1)^{p}$ is spanned by the vector
$
u_{\ka,p+2}+(-1)^{p+1}e^{2\pi ip(p+2)/\ka}u_{\ka,p+6}.
$
So this vector is proportional to $c(\tau)(v_1-v_5)$.  By Lemma \ref{del}, $c(\tau)$ must be a solution of the differential equation 
\begin{multline*}
\frac{\ka}{p+1}\left(\frac{d}{d\tau}c(\tau)\right)(\theta_{6,1}(0,\tau)-\theta_{6,5}(0,\tau))=\\-6\left(\frac{d}{d\tau}\ln\thi'(0,\tau)\right)c(\tau)(\theta_{6,1}(0,\tau)-\theta_{6,5}(0,\tau))+10c(\tau)\frac{d}{d\tau}(\theta_{6,1}(0,\tau)-\theta_{6,5}(0,\tau)). 
\end{multline*}
The function $c(\tau)=\left((2\pi)^{3}(\theta_{6,1}(0,\tau)-\theta_{6,5}(0,\tau))^{5}\thi'(0,\tau)^{-3}\right)^{2(p+1)/\ka}$ is a solution of the preceding equation.
By Lemma \ref{p4},
$c(\tau)=(\eta(\tau))^{-8(p+1)/\ka}$.
The coefficient of proportionality is computed as in the proof of \Ref{esi1}.
Notice that the solution in \Ref{esi10} is invariant with respect to the action of the modular group.\qed 

\end{document}